\documentclass[12pt,twoside]{article}
\usepackage{amsmath}
\usepackage{amsfonts}
\usepackage{amsthm}
\usepackage{amssymb}
\usepackage{graphicx}
\usepackage{multicol}

\begin{document}
\title{{\normalsize
{\bf Classifying the surface-knot modules}}}
\author{{\footnotesize Akio Kawauchi}\\
\date{}
{\footnotesize{\it Osaka Central Advanced Mathematical Institute, Osaka Metropolitan University}}\\ 
{\footnotesize{\it Sugimoto, Sumiyoshi-ku, Osaka 558-8585, Japan}}\\ 
{\footnotesize{\it kawauchi@omu.ac.jp}}}
\maketitle
\vspace{0.25in}
\baselineskip=10pt
\newtheorem{Theorem}{Theorem}[section]
\newtheorem{Conjecture}[Theorem]{Conjecture}
\newtheorem{Lemma}[Theorem]{Lemma}
\newtheorem{Sublemma}[Theorem]{Sublemma}
\newtheorem{Proposition}[Theorem]{Proposition}
\newtheorem{Corollary}[Theorem]{Corollary}
\newtheorem{Claim}[Theorem]{Claim}
\newtheorem{Definition}[Theorem]{Definition}
\newtheorem{Example}[Theorem]{Example}

\begin{abstract} The $k$th module of a surface-knot of a genus $g$ in the 4-sphere is the 
 $k$th integral homology module of the infinite cyclic covering of the surface-knot 
 complement.  The reduced first module is the quotient module of the first module by the finite 
 sub-module  defining the torsion linking. It is shown that the reduced first module for every 
 genus $g$ is  characterized in terms of properties of a finitely generated module. As a 
 by-product,  a concrete example of the fundamental group of a surface-knot of genus $g$ 
 which is not the fundamental group of any surface-knot of genus $g-1$ is given for every 
 $g>0$. The torsion  part and the torsion-free part of the second module are determined by 
 the reduced first module  and the genus-class on the reduced first module. The third module 
 vanishes. The concept of an exact leaf of a surface-knot is introduced, whose linking  
is an orthogonal sum of  the torsion linking and a hyperbolic linking.  

\phantom{x}



\end{abstract}

\baselineskip=15pt

\bigskip

\noindent{\bf 1. Introduction}

A {\it surface-knot} is a closed (connected oriented) surface $F$ with genus $g (\geq 0)$ 
smoothly embedded in the 4-sphere $S^4$. 
Let $E=E(F)=\mbox{cl}(S^4\setminus N(F))$ be the exterior of a surface-knot $F$,  where 
$N(F)=F\times D^2$ is a normal disk bundle of $F$ in $S^4$,  where the section 
$F\times 1$ of the circle bundle $\partial N(F)=\partial E=F\times S^1$ of $F$ is chosen so that the natural homomorphism $H_1(F\times1;Z)\to H_1(E;Z)=Z$ is the zero map. 
Let $\mbox{proj}:\tilde E\to E$ be the infinite cyclic connected covering 
belonging to the kernel of the canonical epimorphism $\pi_1(E, x_0)\to H_1(E;Z)=Z$. Then 
the section $F\times 1$ of $\partial E$ lifts to the section $F\times 0$ of 
$\partial E=F\times\mbox{\boldmath $R$}$. 
Let $\Lambda=Z[Z]=Z[t, t^{-1}]$ be the integral group ring of the infinite cyclic covering transformation group $<t>$ of $\tilde E$ with generator $t$ identified with the meridian generator of $F$ in $H_1(E;Z)=Z$. 
The $k${\it th surface-knot module} (or simply the $k${\it th module} of a surface-knot $F$ in $S^4$ is the $k$th integral homology group $A_k(F)=H_k(\tilde E;Z)$ considered as a finitely generated $\Lambda$-module. 
For a finitely generated $\Lambda$-module $H$,  let $TH$ be the $\Lambda$-torsion part of $H$ and $BH=H/TH$, the $\Lambda$-torsion-free part of $H$. Let $DH$ be the $\Lambda$-submodule of 
$TH$ consisting of every element $x$ with $f_i(t)x=0\,(i=1,2,\dots,s)$ for a coprime 
element system $f_i(t)\in\Lambda\,(i=1,2,\dots,s)$, which is the maximal finite $\Lambda$-submodule of $TH$,  and $T_DH=TH/DH$. Let $E^q(H)=Ext^q_{\Lambda}(H, \Lambda)$ be the $q$th extension 
cohomology $\Lambda$-module of $H$. Since $\Lambda$ is a Noetherian ring of global dimension $2$,  $E^q(H)$ is always finitely generated and $E^q(H) =0\,  (q\geq 3)$. In particular,  
$E^0(H)= \hom_{\Lambda}(H, \Lambda)$ is a free 
$\Lambda$-module,  whose $\Lambda$-rank is defined to be the $\Lambda$-{\it rank} of $H$. It is a standard fact that there is a natural short exact sequence 
\[0\to E^1(BH)\to E^1(H)\to E^1(TH)\to 0, \]
where $E^1(BH)$ is a finite $\Lambda$-module and $E^1(H)$ is a finitely generated torsion 
$\Lambda$-module with 
\[E^1(BH)\cong DE^1(H),\quad T_DE^1(H)=\hom_{\Lambda}(TH,  Q(\Lambda)/\Lambda)
=\hom_{\Lambda}(T_DH,  Q(\Lambda)/\Lambda)\] 
for the quotient field $Q(\Lambda)$ of $\Lambda$ and $E^1E^1(H)=E^1E^1(T_DH)=T_DH$. 
The $\Lambda$-module $E^2(H)$ is a finite $\Lambda$-module 
with $E^2(H)=\hom_Z(DH, Q/Z)$ and $E^2E^2(H)=E^2E^2(DH)=DH$. 
It is also a standard fact that there is a natural short exact sequence 
\[0\to BH\to E^0E^0(BH)\to E^2E^1(BH)\to 0.\]
A $(t-1)$-{\it divisible} $\Lambda$-module is a finitely generated $\Lambda$-module 
$H$ such that the multiplication $t-1:H\to H$ is a $\Lambda$-isomorphism. 
Then every $\Lambda$-submodule and every quotient $\Lambda$-module of $H$ are 
torsion $(t-1)$-divisible $\Lambda$-modules and $DH$ is a finite $\Lambda$-module. 
See \cite{K86, Kerv,  Levine} for these properties of $E^q(H)$. 

An $r$-{\it weight} of a finite $\Lambda$-module $D$ is a $\Lambda$-epimorphism 
$\omega:\Lambda^r\to D$. 
Two $r$-weights $\omega$ and $\omega'$ of $D$ are {\it equivalent} if there are 
$\Lambda$-isomorphisms $f_{\Lambda}:\Lambda^r\to \Lambda^r$ and $f_D: D\to D$ such that $\omega' =f_D\omega f_{\Lambda}^{-1} $. 
An $r$-{\it class} on $D$ is the equivalence class $[\omega]$ of an $r$-weight 
$\omega$ of $D$. For every $r$,  there are only finitely many $r$-classes on $D$,  where if there is no $\Lambda$-epimorphism $\Lambda^r\to D$,  then we understand that $D$ has 
the {\it empty \mbox{$r$}-class} $[\emptyset]$. 
For every non-empty $r$-class $[\omega]$ on $D$, then there is a unique (up to 
$\Lambda$-isomorphisms) torsion-free  $\Lambda$-module  $B$ such that 
the natural $\Lambda$-epimorphism $E^0E^0(B)\to E^2E^1(B)$  is equivalent to $\omega$.
 (see Lemma~4.1). 

Elementary computations show that 
$A_k(F)=0$ except for $0\leq k\leq 3$, and $A_0(F)=Z$ (regarded as a $\Lambda$-module with trivial $t$-action). 
By the zeroth duality of \cite{K86},  
there is a non-degenerate $\Lambda$-Hermitian form 
\[S: E^0E^0(BA_2(F))\times E^0E^0(BA_2(F)) \to \Lambda\]
as an invariant of a surface-knot $F$ in $S^4$ with the identities 
\[f(t)S(x, x')=S(f(t^{-1})x, x')=S(x, f(t)x')\quad (x, x'\in BA_2(F), \,  f(t)\in \Lambda) \] 
extending the non-degenerate $\Lambda$-intersection form
\[S^B:BA_2(F) \times BA_2(F) \to \Lambda \]
defined by 
\[S^B(x, x')=\mbox{Int}_{\Lambda}(x, x')=\sum_{i=-\infty}^{+\infty} 
\mbox{Int}(x, t^i x')t^{-i}\in \Lambda.\] 
By the second duality of \cite{K86},  
the torsion linking (that is a $t$-isometric symmetric bilinear non-singular pairing) 
\[\ell_F :\Theta(F)\times \Theta(F) \to Q/Z\]
on a finite $\Lambda$-module $\Theta(F)$ in $DA_1(F)$ 
is defined as an invariant of a surface-knot $F$ in $S^4$. 
The {\it reduced first module} of $F$ in $S^4$ is 
the quotient $\Lambda$-module $R_1(F)=A_1(F)/\Theta(F)$ of the first module $A_1(F)$,  
which is an invariant of a surface-knot $F$ in $S^4$. 
Let $e(H)$ denote the minimal number of $\Lambda$-generators of $H$.
The following theorem is the main result of this paper.

\phantom{x}

\noindent{\bf Theorem~1.1.} The $k$th surface-knot modules $A_k(F)\,  (1\leq k\leq 3)$ of 
every surface-knot $F$ of genus $g>0$ in $S^4$ have the following properties.

\medskip

\noindent{(1)}
A $\Lambda$-module $H$ is $\Lambda$-isomorphic to the reduced first module $R_1(F)$ of a surface-knot $F$ in $S^4$ of genus $g\,  (\geq 0)$ if and only if $H$ is a $(t-1)$-divisible 
finitely generated $\Lambda$-module with inequality $e(E^2(H))\leq g$.

\medskip

\noindent{(2)} Every surface-knot $F$ in $S^4$ of genus $g$ defines a $g$-class 
invariant $[\omega_F]$ on the finite $\Lambda$-module $E^2(R_1(F))$  so  
that the reduced first module $R_1(F)$ and the $g$-class $[\omega_F]$ determine 
the $\Lambda$-modules $TA_2(F)$ and $BA_2(F)$ up to $\Lambda$-isomorphisms. 
In particular,  there are $t$-anti $\Lambda$-isomorphisms
\[E^1(R_1(F)) \cong TA_2(F) , \quad E^2(R_1(F))\cong E^2E^1(BA_2(F)).\] 

\medskip

\noindent{(3)} There is a direct sum splitting 
$BA_2(F)=X_F\oplus Y_F$ with 
$Y_F$ a free $\Lambda$-module of rank $g$ such that the $\Lambda$-Hermitian form 
\[S:E^0E^0(BA_2(F))\times E^0E^0(BA_2(F)) \to \Lambda\]
is given by 
\[S(x_i, x_j)=S(y_i, y_j)=0,  \quad S(x_i, y_j)=(t-1)\delta_{ij} \quad (i, j=1, 2, \dots,  g)\] 
for a $\Lambda$-basis $x_i,  y_i \,  (i =1, 2, \dots,  g)$ of 
$E^0E^0(BA_2(F))=E^0E^0(X_F)\oplus Y_F$ with $x_i \,  (i =1, 2, \dots,  g)$ a $\Lambda$-basis of $E^0E^0(X_F)$ and $y_i\,  (i=1, 2, \dots, g)$ a $\Lambda$-basis of $Y_F$.

\medskip

\noindent{(4)} $A_3(F)=0$.

\phantom{x}

The $g$-class  $[\omega_F]$ on the finite $\Lambda$-module $E^2(R_1(F))$ 
is called the {\it genus-class invariant} of  a surface-knot $F$.
The weaker inequality $e(E^2(R_1(F)))\leq 2g$ has been earlier obtained and applied 
to surface-knot theory (see \cite[p.192]{K96}). 
If $F$ is an $S^2$-knot $K$ in $S^4$,  then $e(E^2(R_1(F)))=0$,  that is,  $R_1(K)$ is a
$Z$-torsion-free $\Lambda$-module,  which is also the result of Farber-Levine pairing of an $S^2$-knot in $S^4$ (\cite{Farber, Levine}). 
This weaker inequality and the symmetric property of $\Theta(F)$ that $\Theta(F)$ admits a $t$-anti automorphism are applied to know implicitly the properness of the sequence 
\[\mbox{\boldmath $G$}(0)\subset \mbox{\boldmath $G$}(1)\subset\mbox{\boldmath $G$}(2)\subset \dots \subset \mbox{\boldmath $G$}(g) \subset \dots\]
where $\mbox{\boldmath $G$}(g)$ denotes the set of the fundamental groups of 
surface-knots of genus $g$ (see \cite{K07}) and the properness of the sequence 
\[\mbox{\boldmath $A$}(0)\subset \mbox{\boldmath $A$}(1)\subset \mbox{\boldmath $A$}(2)\subset \dots \subset \mbox{\boldmath $A$}(g) \subset \dots\]
where $\mbox{\boldmath $A$}(g)$ denotes the set of the first modules of surface-knots of genus $g$ (see \cite{K08}). 
By Theorem~1.1 (1) and the symmetric property of $\Theta(F)$,  
the properness of these sequences can be shown with explicit examples as follows. 

\phantom{x}

\noindent{\bf Corollary~1.2.} For every prime $p \geq 5$,  consider the finite $\Lambda$-module 
$D=\Lambda/(p,  2t-1)$ and the ribbon presented group
\[\pi =<x,  y |\,  y = (x^{-1}y) x (y^{-1}x),  \,  y=(x y^{-1})^p y (y x^{-1})^{p} >.\]
Then there is a ribbon torus-knot $T$ in $S^4$ with fundamental group 
$\pi_1(S^4\setminus T, x_0)=\pi$ 
and $A_1(T)=D$. For every integer $ g\geq 1$,  let $T_g$ be the $g$-fold connected sum of $T$ in $S^4$ 
which is a ribbon surface-knot of genus $g$. Then the fundamental 
group $\pi_1(S^4\setminus T_g, x_0)$ which has a  ribbon presentation 
\[<x,  y_1,y_2,\dots,y_g |\,  
y_i = (x^{-1}y_i) x (y_i^{-1}x),\,  y_i=(x y_i^{-1})^p y_i (y_i x^{-1})^{p},\, i=1,2,\dots,g>\]
belongs to 
$\mbox{\boldmath $G$}(g)\setminus \mbox{\boldmath $G$}(g-1)$ and 
the first module $A_1(T_g)$ of $T_g$ in $S^4$ belongs to 
$\mbox{\boldmath $A$}(g)\setminus \mbox{\boldmath $A$}(g-1)$. 

\phantom{x}

A basic idea of the proof of Theorem~1.1 is to construct a surface-basis for every 
surface-knot $F$ of genus $g$ in $S^4$ to apply the 3 dualities in \cite{K86} which is described from now.
A {\it loop basis} for a closed oriented surface $F$ of genus $g>0$ is a system of simple loops $ \alpha_i ,  \alpha'_i \,  (i=1,  2, \dots ,  g)$ in $F$ such that 
\[\alpha_i \cap \alpha_j = \alpha_i \cap \alpha'_j = \emptyset\, (i\ne j) \quad\mbox{and} 
\quad \alpha_i \cap \alpha'_i =p_i, \mbox{a point}.\]
A loop basis $\alpha_i ,  \alpha'_i \,  (i=1,  2, \dots ,  g)$ of a surface-knot $F$ is {\it spin} 
if $q([\alpha_i ]_2)=q([\alpha'_i ]_2)=0$ for all $i$ with respect to the quadratic function 
$q: H_1(F;Z_2)\to Z_2$ associated with the surface-knot $F$ in $S^4$.
By \cite{HK},  there is a spin loop basis for every surface-knot $F$ in $S^4$. 

\phantom{x}

\noindent{\bf Definition.}
A {\it surface-basis } of a surface-knot $F$ in $S^4$ of genus $g>0$ is a system 
of (compact connected oriented) surfaces $D_i,  D'_i\,  (i=1, 2, \dots , g)$ smoothly embedded in $S^4$ such that 

\medspace 

\noindent{(1)} $ D_i \cap F=\partial D_i= \alpha_i$ and $D'_i \cap F= \partial D'_i= \alpha'_i , \, (i=1, 2, \dots, g))$ for a spin loop basis $ \alpha_i ,  \alpha'_i \,  (i=1,  2, \dots ,  g)$ of $F$,  

\medspace 

\noindent{(2)} $D_i \cap D_j = D'_i \cap D'_j = D_i \cap D'_j = \emptyset\,  (i\ne j)$,  and
the self $Z$-intersection numbers
$\mbox{Int}(D_i,  D_i)=\mbox{Int}(D'_i, D'_i)=0$ with respect to the surface framing of $F$ 
for all $i$, and

\medspace 

\noindent{(3)} the natural homomorphisms 
$H_1(D_i\setminus \alpha_i;Z)\to H_1(S^4\setminus F;Z)$ 
and $H_1(D'_i\setminus \alpha'_i;Z)\to H_1(S^4\setminus F;Z)$ are the zero maps for all $i$. 

\phantom{x}

In this definition,  note that no information on the intersection between 
the interior $\mbox{Int}D_i$ of $D_i$ and the interior $\mbox{Int}D'_i$ of $D'_i$ is given 
for every $i$. and the interchange between some surfaces in $D_i\,  (i=1, 2, \dots , g)$ and 
the corresponding surfaces in $D'_i\,  (i=1, 2, \dots , g)$ makes a surface-basis for $F$ in $S^4$. 
The following theorem is basically important in this paper,  which is shown in Section~2. 

\phantom{x}

\noindent{\bf Theorem~1.3.} 
For every spin loop system $ \alpha_i ,  \alpha'_i \,  (i=1,  2, \dots ,  g)$ 
of a surface-knot $F$ of genus $g$ in $S^4$,  there is a surface-basis $D_i ,  D'_i\, (i=1, 2, \dots, g)$ for a surface-knot $F$ in $S^4$ with 
$\partial D_i=\alpha_i ,  \partial D'_i=\alpha'_i \,  (i=1,  2, \dots ,  g)$.

\phantom{x}

A {\it leaf} (or in other words,  a {\it Seifert hypersurface}) of a surface-knot $F$ in $S^4$ 
is a compact connected oriented 
3-manifold $V_F$ (smoothly embedded) in $S^4$ with $\partial V_F=F$,  which is always exists (see \cite{Gluck},  \cite[II] {KSS}). 
A leaf $V_F$ is also considered as a proper 3-submanifold of $E$ with $\partial V_F=F\times 1\subset F\times S^1=\partial N(F)$. 
Then the homology class $[V_F]\in H_3(\tilde E, \partial \tilde E;Z)$ is just 
the {\it fundamental class} of the covering $\mbox{proj}: \tilde E\to E$ (see \cite{K98}). 
A leaf $V_F$ of $F$ in $E$ is {\it exact} if the sequence 
\[0\to \mbox{Tor}H_2(\tilde E, \tilde V_F;Z)\to \mbox{Tor}H_1(\tilde V_F;Z) 
\to \mbox{Tor}H_1(\tilde E;Z)\]
is exact. This notion is a variation of a closed exact leaf on a closed oriented 4-manifold with infinite cyclic first homology group in \cite{K00}.

\phantom{x}

\noindent{\bf Theorem~1.4.} For every surface-basis $D_i ,  D'_i\, (i=1, 2, \dots, g)$ of every 
surface-knot $F$ of genus $g$ in $S^4$,  there is an exact leaf $V_F$ containing 
the half surface-basis $D_i \, (i=1, 2, \dots, g)$ as proper surfaces. 

\phantom{x}

A {\it hyperbolic linking} is a linking (i.e.,  non-singular symmetric bilinear form) 
$\ell:G^2\times G^2\to Q/Z$ on the direct double $G^2$ of a finite abelian group $G$ 
such that $\ell(x, x)=0$ for all $x\in G$ (see \cite{KK}). 
The following corollary is a combination result of Theorem~1.4 and an earlier result on 
a closed exact leaf in \cite{K00+}.

\phantom{x}

\noindent{\bf Corollary~1.5.} The torsion linking 
$\ell_F :\Theta(F)\times \Theta(F) \to Q/Z$ 
of every surface-knot $F$ in $S^4$ is an orthogonal summand of the linking 
$\ell_V: \mbox{Tor}H_1(V_F;Z)\times \mbox{Tor}H_1(V_F;Z)\to Q/Z$ for every exact leaf $V_F$ 
containing the half surface-basis $D_i \, (i=1, 2, \dots, g)$ of every surface-basis 
$D_i,  D'_i \, (i=1, 2, \dots, g)$ as proper surfaces,  
which is a non-singular linking and whose complement linking is a hyperbolic linking.

\phantom{x}

In Section~2,  a surface-basis for every surface-knot is constructed.
In Section~3,  the surface-knot manifold $M$ which is a closed spin 4-manifold with 
$H_1(M;Z)\cong Z$ obtained from $S^4$ by a surgery along the surface-knot $F$ is 
considered to apply the 3 dualities of \cite{K86} to the integral infinite cyclic covering homology $H_*(\tilde M;Z)$ where a surface-basis of a surface-knot is used. 
In Section~4,  the proofs of Theorems~1.1 and Corollary~1.2 are given. 
In Section~5,  Theorem~1.4 and Corollary~1.5 are shown by using a closed exact leaf of 
the surface-knot manifold $M$ is discussed in \cite{K00,  K00+}.

\phantom{x}

\noindent{\bf 2. A surface-basis of a surface-knot } 

A {\it surface-basis in the weak sense} for a surface-knot $F$ of genus $g>0$ in $S^4$  is 
a surface-basis for $F$ that does not impose the condition (3).
Namely,  there are (compact connected oriented) surfaces $D_i ,  D'_i\, (i=1, 2, \dots, g)$ smoothly embedded in $S^4$ such that 

\medspace 

\noindent{(1)} $ D_i \cap F=\partial D_i= \alpha_i$ and $D'_i \cap F= \partial D'_i= \alpha'_i , \, (i=1, 2, \dots, g)$ for any given spin loop basis $ \alpha_i ,  \alpha'_i \,  (i=1,  2, \dots ,  g)$ of $F$,  and

\medspace 

\noindent{(2)} the intersection numbers $\mbox{Int}(D_i , D_j )= \mbox{Int}(D'_i ,  D'_j) = 
\mbox{Int}(D_i,  D'_j) =\mbox{Int}(D'_i,  D_j)= 0\,  (i\ne j)$,  and the self $Z$-intersection numbers
$\mbox{Int}(D_i, D_i)=\mbox{Int}(D'_i, D'_i)=0$ with respect to the surface framing of $F$ are $0$ for all $i$ in $S^4$.

\phantom{x}

A surface-basis in the weak sense is constructed  in \cite{HK} for every surface-knot $F$ in $S^4$ with any given spin loop basis $\alpha_i ,  \alpha'_i \,  (i=1,  2, \dots,  g)$.   
To be precise,  the condition that $\mbox{Int}(D_i,  D'_j) = 0\,  (i\ne j)$ is omitted in \cite{HK},  but it is shown as well. 
For a spin loop basis $\alpha_i ,  \alpha'_i \,  (i=1,  2, \dots ,  g)$ in 
$F\times 1\subset\partial E$, let  $D_i ,  D'_i\, (i=1, 2, \dots, g)$ be a surface-basis in the weak sense in $E$ with $\partial D_i= \alpha_i,\, \partial D'_i= \alpha'_i , \, (i=1, 2, \dots, g)$. 
Let $T_i=\ell_i\times S^1,  T'_i=\ell'_i\times S^1\,  (i=1,  2, \dots,  g)$ be the tori in $\partial E$. Let $a(D_i)=T_i\cup D_i$ and $a(D'_i)=T'_i\cup D_i$ be 
be the 2-cycles in $E$ homologous to $T_i$ and $T'_i$,  respectively. 
The elements $[a(D_i)],  [a(D'_i)]\,  (i=1, 2, \dots, g)$ form a basis of $H_2(E;Z)$ whose dual basis of $H_2(E, \partial E;Z)$ with respect to the non-singular  intersection form
\[\mbox{Int}_{\partial}: H_2(E;Z)\times H_2(E, \partial E;Z)\to Z\]
are given by the homology classes $[D'_i],  [D_i]\,  (i=1, 2, \dots, g)$. 
For the homological argument on the infinite cyclic covering $\tilde E$ of 
the exterior $E$ of a surface-knot $F$ of genus $g$,  the following facts will be used throughout the paper:

\phantom{x}

\noindent{\bf Lemma~2.1.} The exterior $E$ of a surface-knot $F$ of genus $g$ 
has the following homological properties. 

\noindent{(1)}  $A_1(F)$ and $H_1(\tilde E, \partial\tilde E;Z)$
are $(t-1)$-divisible and there is a natural 
$\Lambda$-isomorphism $A_1(F)\cong H_1(\tilde E, \partial\tilde E;Z)$. 

\medskip

\noindent{(2)} $DA_2(F)=DH_2(\tilde E, \partial\tilde E;Z)=0$.  

\medskip

\noindent{(3)}
$TA_2(F)$ and $TH_2(\tilde E, \partial\tilde E;Z)$ 
are $(t-1)$-divisible,  so that there is a natural $\Lambda$-isomorphism 
$TA_2(F)\cong TH_2(\tilde E, \partial\tilde E;Z)$ and there is a natural short 
exact sequence 
\[0\to BA_2(F)\stackrel {j_*}{\to} BH_2(\tilde E, \partial\tilde E;Z)\stackrel {\partial_*}{\to} H_1(\partial\tilde E;Z)\to 0\]
with $H_1(\partial\tilde E;Z)\cong Z^{2g}$. 

\medskip

\noindent{(4)}
$E^1(BA_2(F))$ and $E^1(BH_2(\tilde E, \partial\tilde E;Z))$ are $(t-1)$-divisible and there is a natural $\Lambda$-isomorphism 
$E^1(BH_2(\tilde E, \partial\tilde E;Z))\cong E^1(BA_2(F))$. 

\phantom{x}

Technically, the following observation is useful  (whose proof is direct and omitted).

\phantom{x}

\noindent{\bf Observation~2.2.} In an exact sequence 
\[H_0 \to H_1 \stackrel {\varphi}{\to} H_2\to H_3\]
of finitely generated $\Lambda$-modules $H_i\,(0\leq i\leq 3)$ and 
$\Lambda$-homomorphisms, if $(t-1)H_0=(t-1)H_3=0$ and $H_1, H_2$ are  
$(t-1)$-divisible, then  the $\Lambda$-homomorphism $\varphi:H_1\to H_2$ is a 
$\Lambda$-isomorphism.

\phantom{x}

\noindent{\bf Proof of Lemma~2.1.}
Since $H_1(E;Z)\cong Z$,  the Wang exact sequence shows that $t-1:A_1(F)\to A_1(F)$ is an isomorphism,  showing that $A_1(F)$ is 
$(t-1)$-divisible. This fact and the homology exact sequence of the pair 
$(\tilde E, \tilde\partial E)$ shows that $TA_2(F)\cong TH_2(\tilde E, \partial\tilde E;Z)$,  showing (1). By the second duality of \cite{K86},  
there are $t$-anti epimorphisms 
\[\theta: DA_2(F)\to E^1(BH_1(\tilde E, \partial \tilde E;Z))=0,\quad  
\theta: DH_2(\tilde E, \partial \tilde E;Z)\to E^1(BA_1(F))=0\]
whose kernels 
$DA_2(F)^{\theta}=DA_2(F),\,
DH_2(\tilde E, \partial \tilde E;Z)^{\theta}=DH_2(\tilde E, \partial \tilde E;Z)$ 
are $t$-anti $\Lambda$-isomorphic to 
$\hom_Z(DH_0(\tilde E, \partial \tilde E;Z)^{\theta}, Q/Z)=0$ 
and $\hom_Z(DA_0(F)^{\theta}, Q/Z)=0$,  
for $DA_0(F)=DH_0(\tilde E, \partial \tilde E;Z)=0$. 
Hence $DA_2(F)=DH_2(\tilde E, \partial \tilde E;Z)=0$, showing (2). Then 
by the second duality of \cite{K86},  
$TA_2(F)$and $TH_2(\tilde E, \partial \tilde E;Z)$ are $t$-anti $\Lambda$-isomorphic to 
$E^1(T_DH_1(\tilde E, \partial \tilde E;Z))=E^1(H_1(\tilde E, \partial \tilde E;Z))$ and $E^1(T_DA_1(F))=E^1(A_1(F))$ which are $(t-1)$-divisible,  respectively. 
This $(t-1)$-divisibility and the homology exact sequence of the pair 
$(\tilde E, \partial\tilde E)$ show that 
the natural homomorphism 
$TA_2(F)\to TH_2(\tilde E, \partial\tilde E;Z)$ is a 
$\Lambda$-isomorphism,  so that there is a natural short 
exact sequence 
\[(*)\qquad 0\to BA_2(F)\stackrel {j_*}{\to} BH_2(\tilde E, \partial\tilde E;Z)\stackrel {\partial_*}{\to} H_1(\partial\tilde E;Z)\to 0,\]
showing (3). 
By the second duality of \cite{K86},  there are  $t$-anti $\Lambda$-eimorphisms 
\[\theta:DA_1(F)\to E^1(BH_2(\tilde E, \partial\tilde E;Z)),\, 
\theta: DH_1(\tilde E, \partial\tilde E;Z)\to E^1(BA_2(F)).\]
Since $DA_1(F), DH_1(\tilde E, \partial\tilde E;Z)$ are $(t-1)$-divisible by (1), 
it is seen that a natural $\Lambda$-homomorphism  
$E^1(BH_2(\tilde E, \partial\tilde E;Z))\to E^1(BA_2(F))$ 
is a $\Lambda$-isomorphism by applying the extension cohomology to the short exact sequence $(*)$, showing (4). 
This completes the proof of Lemma~2.1.

\phantom{x}

By the zeroth duality of \cite{K86},  the non-singular $\Lambda$-form
\[S_{\partial}: E^0E^0(BA_2(F))\times E^0E^0(BH_2(\tilde E,  \partial \tilde E;Z))\to\Lambda\] 
is given by extending the non-degenerate $\Lambda$-intersection form
the non-degenerate $\Lambda$-Hermitian form
\[S^B_{\partial}:BA_2(F) \times BH_2(\tilde E,  \partial \tilde E;Z) \to \Lambda,  \]
defined by 
$S^B_{\partial}(x, x')=\mbox{Int}_{\Lambda}(x, x')=\sum_{i=-\infty}^{+\infty} 
\mbox{Int}(x, t^i x')t^{-i}\in \Lambda$  
for $x\in BA_2(F)$,  $x'\in BH_2(\tilde E,  \partial \tilde E;Z)$.

Let $\Lambda^+$ be the subring of the quotient field $Q(\Lambda)$ of $\Lambda$ 
generated by the products $u(t)^{-1}f(t)$ of any elements 
$u(t), f(t)\in\Lambda$ with $u(1)=\pm1$. 
Note that the ring $\Lambda^+$ admits a $t$-anti automorphism. 
For a finitely generated $\Lambda$-module $H$,  let 
$H^+=H\otimes_{\Lambda}\Lambda^+$. 
It is a standard fact that {\it for every \mbox{$(t-1)$}-divisible finitely generated 
\mbox{$\Lambda$}-module \mbox{$H$},  there is an element 
\mbox{$u(t)\in\Lambda$} such that \mbox{$u(1)=\pm1$} and  
\mbox{$u(t)H=0$}.} (In fact, $H=TH$. Since $T_DH$ has projective dimension $1$, there 	is a short exact sequence 
$0\to \Lambda^m\stackrel {P(t)}{\to}  \Lambda^m\to T_DH\to 0$, where $P(t)$ denotes 
a presentation $\Lambda$-matrix.
Then $u_1(t)=\mbox{det}P(t)$ has $u_1(t)T_DH=0$. Since $T_DH$ is (t-1)-divisible, 
$P(1)$ is a unimodular matrix and $u_1(1)=\mbox{det}P(1)=\pm1$. 
On the other hand, some iteration 
$(t-1)^{m'}$ of  $t-1$  acts  on the finite $\Lambda$-module $DH$ as the identity. 
Then $u_2(t)=1-(t-1)^{m'}$ has $u_2(1)=1$  and $u_2(t)DH=0$. 
The product $u(t)=u_1(t)u_2(t)$ has $u(1)=\pm1$ and $u(t)H=0$, as desired.)
Since $A_1(F)\cong H_1(\tilde E, \partial \tilde E;Z)$ is $(t-1)$-divisible,  the second duality of \cite{K86} implies that $E^2E^1(BA_2(F))$ and 
$E^2E^1(BH_2(\tilde E,  \partial \tilde E;Z))$ are $(t-1)$-divisible,  so that 
$BA_2(F)^+=E^0E^0(BA_2(F))^+$ and 
$BH_2(\tilde E,  \partial \tilde E;Z)^+=E^0E^0(BH_2(\tilde E,  \partial \tilde E;Z))^+$ 
are free $\Lambda^+$-modules,  and the non-degenerate $\Lambda$-Hermitian form 
$S^B_{\partial}$ induces a {\it non-singular} $\Lambda^+$-Hermitian form 
\[S^+_{\partial}=\mbox{Int}_{\Lambda^+}: 
BA_2(F)^+ \times BH_2(\tilde E,  \partial \tilde E;Z)^+ \to \Lambda^+\] 
by defining 
\[\mbox{Int}_{\Lambda^+}(x,x')= u(t^{-1})^{-1}u'(t)^{-1}\mbox{Int}_{\Lambda}(u(t)x,u'(t) x')\] 
for $x\in BA_2(F)^+, x'\in BH_2(\tilde E, \partial \tilde E;Z)^+$ and $u(t) ,u'(t)\in\Lambda$ 
such that $u(1)=u'(1)=1$ and $u(t)x\in BA_2(F), u'(t)x'\in BH_2(\tilde E,  \partial \tilde E;Z)$. 
Similarly,  the non-degenerate $\Lambda$-Hermitian form 
$S^B: BA_2(F)\times BA_2(F)\to \Lambda $ induces a non-degenerate
$\Lambda^+$-Hermitian form 
\[S^+=\mbox{Int}_{\Lambda^+}: BA_2(F)^+ \times BA_2(F)^+ \to \Lambda^+.\] 
Note that there is a natural short exact sequence 
\[0\to BA_2(F)^+ \stackrel {i_*}{\to} 
BH_2(\tilde E,  \partial \tilde E;Z)^+ \stackrel {\partial_*}{\to} H_1(\partial \tilde E;Z) \to 0\]
and $H_1(\partial \tilde E;Z)=Z^{2g}$ with the $Z$-basis represented 
by the spin loop basis $\alpha_i\times 0 ,  \alpha'_i\times 0 \,  (i=1,  2, \dots ,  g)$ of 
$F\times 0\subset F\times \mbox{\boldmath $ R$} =\partial \tilde E$. 
Note that 
\[S^+(x, x')=S^+_{\partial}(x, i_*(x'))\]
for all $x, x'\in BA_2(F)^+$. 
A {\it well-defined pair} of relative $2$-cycles in $BH_2(\tilde E,  \partial \tilde E;Z)^+$ 
is a pair $(c, c')$ of relative $2$-cycles $c, c'$ in $BH_2(\tilde E,  \partial \tilde E;Z)^+$ 
such that the boundary 1-cycle pair $(\partial c,  \partial c')$ is any pair of 
$\pm\alpha_i\times 0,  \pm\alpha'_i\times 0 \,  (i=1,  2, \dots ,  g)$ except for 
the unordered pair of 
$\pm\alpha_i\times 0$ and $\pm\alpha'_i\times 0$ for every $i$. 
For every well-defined pair $(c, c')$,  the $\Lambda^+$-intersection number 
$\mbox{Int}_{\Lambda^+}(c, c')\in \Lambda^+$ is well-defined where 
$\mbox{Int}_{\Lambda^+}(c, c')$ with $\partial c=\pm \partial c'$ is understood as 
the $\Lambda^+$-intersection number by using by the surface-framing in $F\times 0$. 
Then the following identities hold.
\[(t^{-1}-1) \mbox{Int}_{\Lambda^+}(c, c') =
\mbox{Int}_{\Lambda^+}((t-1)c, c')=S^+_{\partial}(i_*^{-1}[(t-1)c], [c']),\]  
\[(t^{-1}-1)(t-1)\mbox{Int}_{\Lambda^+}(c, c')=S^+(i_*^{-1}[(t-1)c], i_*^{-1}[(t-1)c']).\]
The following lemma is used for the present argument.

\phantom{x}

\noindent{\bf Lemma~2.2.} 
Let $C:F\times[0, 1]\to S^4\times[0, 1]$ be a smooth concordance form 
a surface-knot $F=C(F\times 0)$ with a spin loop system 
$\alpha_i ,  \alpha'_i \,  (i=1,  2, \dots ,  g)$ to a surface-knot 
$G=C(F\times 1)$ in $S^4$ with a spin loop system 
$\beta_i ,  \beta'_i \,  (i=1,  2, \dots ,  g)$. 
Then there is a $\Lambda^+$-isomorphism $\varphi$ from 
the non-singular $\Lambda^+$-Hermitian form 
$S^+_{\partial}: BA_2(F)^+ \times BH_2(\tilde E(F),  \partial \tilde E(F);Z)^+ \to \Lambda^+$ 
to the non-singular $\Lambda^+$-Hermitian form 
$S^+_{\partial}:BA_2(G)^+ \times BH_2(\tilde E(G),  \partial \tilde E(G);Z)^+ \to \Lambda^+$ 
sending the homology classes $[\alpha_i\times 0] ,  [\alpha'_i\times 0]\,  (i=1,  2, \dots ,  g)$ in $H_1(\partial \tilde E(F);Z))$ to the homology classes 
$[\beta_i\times 0] ,  [\beta'_i\times 0]\,  (i=1,  2, \dots ,  g)$ in $H_1(\partial \tilde E(G);Z))$,  
respectively.

\phantom{x}

\noindent{\bf Proof of Lemma~2.2.} 
Let $E(C)=\mbox{cl}(S^4\times[0, 1]\setminus N(F\times[0, 1]))$ be the exterior of 
the concordance $C$. Then $(E(C);E(F), E(G))$ is a homology cobordism with 
$(\partial' E(C); \partial E(F),  \partial E(G))$ the product cobordism for 
$\partial' E(C)=\mbox{cl}(\partial E(C)\setminus (E(F)\cup E(G))$. 
Then $H_*(\tilde E(C),  \tilde E(F);Z)$ and $H_*(\tilde E(C),  \tilde E(G);Z)$ are 
$(t-1)$-divisible 
finitely generated $\Lambda$-modules. Hence 
\[H_*(\tilde E(C),  \tilde E(F);Z)^+=H_*(\tilde E(C),  \tilde E(G);Z)^+=0.\]
Then an argument of the $\Lambda^+$-homology cobordism
 $(\tilde E(C);\tilde E(F),  \tilde E(G))$ similar to the standard homology cobordism argument 
shows the desired result. This completes the proof of Lemma~2.2. ]

\phantom{x}

The proof of Theorem~1.3 is given as follows.

\phantom{x}

\noindent{\bf 2.3: Proof of Theorem~1.3.} 
{\it Every surface-knot \mbox{$F$} in \mbox{$S^4$} is concordant to a trivial surface-knot \mbox{$G$} in \mbox{$S^4$} by a concordance sending any given spin loop basis 
\mbox{$\alpha_i ,  \alpha'_i\,  (i=1,  2, \dots, g)$} of \mbox{$F$} to the standard spin loop basis \mbox{$\beta_i ,  \beta'_i \,  (i=1,  2, \dots,  g)$} of \mbox{$G$}. }
To see this, consider a trivial surface-knot $\bar G$ obtained from $F$ by adding 1-handles 
$h_j\,(j=1,2,\dots,m)$ (see \cite{HoK}). Let $\alpha_i ,  \alpha'_i\,  (i=1,  2, \dots, g)$, 
$\gamma_j ,  \gamma'_j\,  (j=1,  2, \dots, m)$ be a spin loop basis of $\bar G$ with 
$\gamma_j$ a belt loop of $h_j$. By \cite{Hiro}, there is an orientation-preserving diffeomorphism  $f$ of  $(S^4,\bar G)$ sending  the spin loop basis 
$\alpha_i ,  \alpha'_i\,  (i=1,  2, \dots, g)$, $\gamma_j, \gamma'_j$ $(j=1, 2,\dots,m)$
to a standard spin loop basis of $\bar G$.  Thus, there are 2-handles on $\bar G$ 
attached along the loops $\gamma'_j\,  (j=1,  2, \dots, m)$ to obtain  a trivial surface-knot $G$ by the  surgery. Then  $G$ has a spin loop basis $\beta_i ,  \beta'_i \,  (i=1,  2, \dots,  g)$ inherited from the spin loop basis $\alpha_i ,  \alpha'_i\,  (i=1,  2, \dots, g)$. (This is a similar consideration to \cite[(2.5.1), (2.5.1)]{K21}.) The surgery trace gives a desired concordance. 
Let $\Delta_i,  \Delta'_i\,  (i=1,  2, \dots ,  g)$ be a standard disk-basis of $G$ with 
$\partial \Delta_i=\beta_i\times 1,  \partial \Delta'_i=\beta'_i\times 1\,  (i=1,  2, \dots ,  g)$ 
in $G\times1\subset G\times S^1=\partial E(G)$. 
Let $\tilde \Delta_i,  \tilde \Delta'_i\,  (i=1,  2, \dots ,  g)$ be the connected lifts 
of $\Delta_i,  \Delta'_i\,  (i=1,  2, \dots ,  g)$ to $\tilde E(G)$ with 
$\partial \tilde\Delta_i=\beta_i\times 0,  \partial \tilde\Delta'_i=\beta'_i\times 0\,  (i=1,  2, \dots ,  g)$ in $G\times 0\subset G\times \mbox{\boldmath {$R$}}=\partial \tilde E(G)$.
By Lemma~2.2,  there are relative 2-cycles $c_i,  c'_i\,  (i=1,  2, \dots ,  g)$ in 
$BH_2(\tilde E(F),  \partial \tilde E(F);Z)^+$ 
with $\partial c_i=\alpha_i\times 0,  \partial c'_i=\alpha'_i\times 0\,  (i=1,  2, \dots ,  g)$
such that the homology classes $[c_i],  [c'_i]\,  (i=1,  2, \dots ,  g)$ are sent to the 
homology classes $[\tilde\Delta_i],  [\tilde\Delta'_i]\,  (i=1,  2, \dots ,  g)$ in $BH_2(\tilde E(G),  \partial \tilde E(G);Z)^+$ by the $\Lambda^+$-isomorphism $\varphi$. 
Since any pair of $c_i,  c'_i\,  (i=1,  2, \dots ,  g)$ except for $(c_i ,  c'_i),  (c'_i, c_i),  (i=1, 2, \dots, g)$ 
is a well-defined pair,  the following identities hold.
\[\mbox{Int}_{\Lambda^+}(c_i,  c_j)=\mbox{Int}_{\Lambda^+}(\tilde\Delta_i,  \tilde\Delta_j)=0,  \quad
\mbox{Int}_{\Lambda^+}(c'_i,  c'_j)=\mbox{Int}_{\Lambda^+}(\tilde\Delta'_i,  \tilde\Delta'_j)=0\]
for all $i, j$ and
\[\mbox{Int}_{\Lambda^+}(c_i,  c'_j)=\mbox{Int}_{\Lambda^+}(\tilde\Delta_i,  \tilde\Delta'_j)=0, \quad
\mbox{Int}_{\Lambda^+}(c'_i,  c_j)=\mbox{Int}_{\Lambda^+}(\tilde\Delta'_i,  \tilde\Delta_j)=0\]
for all $i, j$ with $i\ne j$. 
There is an element $u(t)\in\Lambda$ with $u(1)=1$ such that the products 
$u(t) [c_i],  u(t) [c'_i]\,  (i=1,  2, \dots ,  g)$ are in $BH_2(\tilde E(F),  \partial \tilde E(F);Z)$ 
for all $i$. 
Since $u(t)$ acts on $Z^{2g}$ as the identity $u(1)=1$, 
there are compact connected oriented proper smoothly embedded surfaces 
$\tilde D_i,  \tilde D'_i \, (i=1, 2, \dots, g)$ in $\tilde E$ with 
$\partial \tilde D_i=\alpha_i,  \partial \tilde D'_i=\alpha'_i\times 0 \, (i=1, 2, \dots, g)$ in
$F\times 0\subset \partial\tilde E$ such that 
\[u(t) [c_i]=[\tilde D_i],  \quad u(t) [c'_i]=[\tilde D'_i]\quad (i=1, 2, \dots, g)\]
in $BH_2(\tilde E,  \partial \tilde E(F);Z)$. 
The $\Lambda$-intersection numbers 
\[\mbox{Int}_{\Lambda}(\tilde D_i, \tilde D_j)=
\mbox{Int}_{\Lambda}(\tilde D'_i, \tilde D'_j)=0\]
for all $i, j$ and 
\[\mbox{Int}_{\Lambda}(\tilde D_i, \tilde D'_j)=\mbox{Int}_{\Lambda}(\tilde D'_i, \tilde D_j)=0\] 
for every $i, j$ with $i\ne j$. 
Then the proper surfaces $\tilde D_i, \tilde D'_i\, (i=1, 2, \dots, g)$ in $\tilde E$ are modified without 
changing the boundary loops into higher genus surfaces which are embeddable into $E$ under the covering projection $\mbox{proj}: \tilde E\to E$ by \cite[Theorem~4.1]{K00}. 
By writing $\mbox{proj}(\tilde D_i) ,  \mbox{proj}(\tilde D'_i)\, (i=1, 2, \dots, g)$ as 
$D_i ,  D'_i\, (i=1, 2, \dots, g)$,  
a surface-basis $D_i ,  D'_i\, (i=1, 2, \dots, g)$ of $F$ in $S^4$ is obtained. 
This completes the proof of Theorem~1.3. 

\phantom{x}

\noindent{\bf 3. The infinite cyclic covering homology of the surface-knot manifold}

Let $D_i ,  D'_i\, (i=1, 2, \dots, g)$ be a surface-basis of $F$ in $E$ with 
$\partial D_i=\alpha_i\times 1 ,  \partial D'_i=\alpha'_i\times 1\, (i=1, 2, \dots, g)$ in 
$F\times 1\subset F\times S^1=\partial E$ by Theorem~1.3. 
Let $V_0$ be a handlebody of genus $g$ with $\partial V_0=F$ 
such that $\alpha_i \,  (i=1,  2, \dots ,  g)$ bound disjoint disks in $V_0$. 
The {\it surface-knot manifold} of a surface-knot $F$ in $S^4$ is the 4D manifold 
$M=E\cup V_0\times S^1$ obtained from $S^4$ by replacing $N(F)=F\times D^2$ with $V_0\times S^1$. Then $H_1(M;Z)\cong Z$.
Let $a(D'_i)=T_i\cup D'_i\, (i=1, 2, \dots, g)$ be the 2-cycles in $M$ homologous to 
$T_i\, (i=1, 2, \dots, g)$,  and 
$s(D_i)=D_i\cup d_i\times 1\, (i=1, 2, \dots, g)$ the closed connected oriented surfaces for disjoint disks $d_i$ in $V_0\times 1$ with $\partial d_i=\alpha_i\, (i=1, 2, \dots, g)$. 
The second homology $H_2(M;Z)$ is a free abelian group of rank $2g$ with a basis consisting of the homology classes $[a(D'_i)],  [s(D_i)]\, (i=1, 2, \dots, g)$ with intersection numbers 
\[\mbox{Int}([a(D'_i)], [a(D'_j)])=\mbox{Int}([s(D_i)], [s(D_j)])=0, \]
\[\mbox{Int}([a(D_i)], [s(D_j)])=\mbox{Int}([s(D_i)], [a(D_j)])=\delta_{ij}\]
for all $i, j$.

Let $\mbox{proj}: \tilde M\to M$ be the infinite cyclic covering of $M$ with 
$\tilde M=\tilde E\cup V_0\times  \mbox{\boldmath $R$}$. 
Let $\tilde D_i ,  \tilde D'_i\, (i=1, 2, \dots, g)$ be the connected lifts of $D_i ,  D'_i\, (i=1, 2, \dots, g)$
with $\partial \tilde D_i=\alpha_i\times 0 ,  \partial \tilde D'_i=\alpha'_i\times 0\, (i=1, 2, \dots, g)$ 
in $F\times 0\subset F\times  \mbox{\boldmath $R$}=\partial \tilde E$. 
Let 
\[ a(\tilde D'_i)=(-\tilde D'_i)\cup [0, 1]\cup t \tilde D'_i, \quad 
s(\tilde D_i)=\tilde D_i\cup d_i\times 0\qquad (i=1, 2, \dots, g)\]
be the closed connected oriented surfaces in $\tilde M$. 
Let 
\[x'_i=[a(\tilde D'_i)], \quad y_i=[s(\tilde D_i)]\qquad (i=1, 2, \dots, g)\]
be the homology classes in $H_2(\tilde M;Z)$. 
Let $ X$ be the $\Lambda$-submodule of $BH_2(\tilde M;Z)$ generated over $\Lambda$ 
by the elements $ x\in BH_2(\tilde M;Z)$ with $S^B(x, x'_i)=0$ for all $i$,  and $Y$ 
the $\Lambda$-submodule of $BH_2(\tilde M;Z)$ generated over $\Lambda$ 
by the elements $ y_i\,  (i=1, 2, \dots, g)$. The following lemma is shown.

\phantom{x}

\noindent{\bf Lemma~3.1.} There is a direct sum splitting $E^0E^0(X)\oplus Y$ 
of the free $\Lambda$-module $E^0E^0(BH_2(\tilde M;Z))$
with $y_i\,  (i=1, 2, \dots, g)$ a $\Lambda$-basis of $Y$ 
such that the $\Lambda$-Hermitian form 
\[S:E^0E^0(BH_2(\tilde M;Z))\times E^0E^0(BH_2(\tilde M;Z)) \to \Lambda\]
is given by 
\[S(x_i, x_j)=S(y_i, y_j)=0,  \quad S(x_i, y_j)=\delta_{ij} \quad (i, j=1, 2, \dots,  g)\]
for a $\Lambda$-basis $x_i\,  (i=1, 2, \dots, g)$ of $E^0E^0(X)$. 

\phantom{x}

\noindent{\bf Proof of Lemma~3.1.} 
By construction,  $S(x'_i, x'_j)=S(y_i, y_j)=S(x'_i, y_j)=0\,  (i\ne j)$ and $S(x'_i, y_i)=1 +(t-1)f_i(t)$ for some $f_i(t)\in\Lambda$. 
Let $X_i$ be the quotient rank one $\Lambda$-module of $X$ by 
the maximal submodule generated over $\Lambda$ by $x'_j$ for all $j\ne i$,  so that 
$X_i$ is a torsion-free $\Lambda$-module of rank one and $E^0(X_i)\cong \Lambda$. 
Let $q_i\in E^0(X_i)$ be the $\Lambda$-homomorphism sending $x\in X_i$ to $S(x, y_i)\in \Lambda$. 
Then it is shown that the element $q_i$ is a generator of $E^0(X_i)\cong\Lambda$. 
To see this,  under the identification $E^0(X_i)=\Lambda$,  
suppose $q_i$ is a non-unit element $q_i=q_i(t)$ in $\Lambda$. 
Then $q_i(1)=\pm 1$ since $q_i(t)$ divides the polynomial $1 +(t-1)f_i(t)$. 
Let $p$ be a prime number such that $q_i(t)$ is still a non-unit polynomial in 
the principal ideal domain $\Lambda_p=Z_p[t, t^-1]$ and the first $Z$-torsion product 
$\mbox{Tor}_1(H_1(\tilde M;Z), Z_p)=0$ by using that the $Z$-torsion $\Lambda$-submodule 
of $H_1(\tilde M;Z)$ is finite because $H_1(\tilde M;Z)$ is $(t-1)$-indivisible. 
Then the universal coefficient theorem means 
$H_2(\tilde M;Z_p)=H_2(\tilde M;Z)\otimes Z_p$. Hence $X\otimes Z_p$ is a self-orthogonal 
complement with respect to the nonsingular $\Lambda_p$-intersection form 
\[S_p: BH_2(\tilde M;Z_p)\times BH_2(\tilde M;Z_p) \to \Lambda_p\]
in \cite{K77}. This means that there is an element $x''_i$ in $X_i$ such that 
$S(x''_i, y_i)=1+pg_i(t)$ for some element $g_i(t)\in \Lambda$,  so that $q_i(t)$ must be a unit element in $\Lambda_p$,  which contradicts that $q_i(t)$ is a non-unit element in $\Lambda_p$. 
Thus,  $q_i$ is a unit element in $\Lambda$. 
Let $\bar q_i\in E^0(X)$ be the image of $q_i$ under the natural monomorphism 
$E^0(X_i)\to E^0(X)$. Then 
the elements $\bar q_i\, (i=1, 2, \dots, g)$ form a $\Lambda$-basis for $E^0(X)$. In fact,  for 
every element $q\in E^0(X)$,  let $q(x'_i)=c_i(t)$ be the element of $\Lambda$. 
Then $q=\sum_{i=1}^g c_i(t) \bar q_i$. If $\sum_{i=1}^g c'_i(t) \bar q_i=0$,  then 
$c'_i(t)\bar q_i(x'_i)=c'_i(t)(1 +(t-1)f_i(t))=0$ and $c'_i=0$ for all $i$.
Let $\bar q_i^*\in E^0E^0(X) \, (i=1, 2, \dots, g)$ be the dual basis of $\bar q_i \, (i=1, 2, \dots, g)$of $E^0(X)$. Since $S(\bar q_i^*, \bar q_j)=S(y_i, y_j)=0$ and $S( \bar q_i^*, y_j)=\delta_{ij}$ for all $i, j$,  the elements $x_i=\bar q_i^*,  y_i\, (i=1, 2, \dots, g)$ form a desired $\Lambda$-basis for $E^0E^0(BH_2(\tilde M;Z))=E^0E^0(X)\oplus Y$. 
This completes the proof of Theorem~3.1.

\phantom{x}

The following corollary is obtained from the proof of Lemma~3.1.

\phantom{x}

\noindent{\bf Corollary~3.2.} For the elements $x'_i=[a(\tilde D'_i)],  
y_i=[s(\tilde D_i)]\, (i=1, 2, \dots, g)$ in $X\oplus Y=BH_2(\tilde M;Z)$,  an element 
$x\in BH_2(\tilde M;Z)$ belongs to the direct summand $X$ if and only if 
the product $u(t)x$ for an element $u(t)\in\Lambda$ with $u(1)=\pm1$ is a linear combination of $x'_i\, (i=1, 2, \dots, g)$ with coefficients in $\Lambda$.

\phantom{x}

\noindent{\bf Proof of Corollary~3.2.} In the proof of Lemma~2.2,  the identities 
$(1+(t-1)f_i(t))x_i=x_i'\, (i=1, 2, \dots, g)$ hold,  so that if $x\in BH_2(\tilde M;Z)$ is in $X$,  then
the product $u(t)x$ for some $u(t)$ with $u(1)=\pm1$ is a linear combination of $x'_i\, (i=1, 2, \dots, g)$ with coefficients in $\Lambda$. 
Conversely,  since $X$ is self-orthogonal with respect to the non-degenerate 
$\Lambda$-intersection form $S_M: BH_2(\tilde M;Z) \times BH_2(\tilde M;Z) \to \Lambda$ and every linear combination of $x'_i\, (i=1, 2, \dots, g)$ 
with coefficients in $\Lambda$ is in $X$,  if $u(t)x$ for some $u(t)\in\Lambda$ with $u(1)=\pm1$ is in $X$,  then $x$ is in $X$. 
This completes the proof of Corollary~3.2.

\phantom{x}

\noindent{\bf 4. Proofs of Theorems~1.1 and Corollary~1.2}

The following lemma is a classification of finitely generated torsion-free $\Lambda$-modules. 

\phantom{x}

\noindent{\bf Lemma~4.1.} For a finite $\Lambda$-module $D$,  let 
$[\omega_i] \,  (i=1, 2, \dots,  n_r)$ be all the $r$-classes on $D$,  and 
$B_i$ a finitely generated torsion-free $\Lambda$-module of rank $r$ 
given by the kernel $\mbox{ker}(\omega_i:\Lambda^r\to D)$ for every $i$.  
Then $B_i\,  (i=1, 2, \dots,  n_r)$ are mutually distinct up to $\Lambda$-isomorphisms  
and every finitely generated torsion-free $\Lambda$-module $B$ of rank $r$ with 
$E^2E^1(B)\cong D$ is $\Lambda$-isomorphic to $B_i$ for some $i$. 
Further,  any two finitely generated torsion-free 
$\Lambda$-modules $B$ and $B'$ with $E^2E^1(B) \cong E^2E^1(B')$ are stably 
$\Lambda$-isomorphic,  i.e.,  $B\oplus\Lambda^m\cong B'\oplus \Lambda^{m'}$ 
for some non-negative integers $m,  m'$. 

\phantom{x}

\noindent{\bf Proof of Lemma~4.1.} For a finitely generated torsion-free $\Lambda$-module $B$ with $E^2E^1(B)\cong D$,  there is a short exact sequence 
$ 0\to B\to E^0E^0(B)\to E^2E^1(B)\to 0$. Since there are $\Lambda$-isomorphisms 
$g_B:E^0E^0(B)\to \Lambda^r$ and $g_D:E^2E^1(B)\to D$ to define 
an $r$-weight $\omega_B:\Lambda^r\to D$ whose kernel 
$\mbox{Ker}(\omega_B)=B_D$ 
is $\Lambda$-isomorphic to $B$,
If there is a $\Lambda$-isomorphism $f:B\to B'$,  
then the $\Lambda$-isomorphism $f$ induces a $\Lambda$-isomorphism from the short 
exact sequence $0\to B\to E^0E^0(B)\to E^2E^1(B)\to 0$ to the short 
exact sequence $ 0\to B'\to E^0E^0(B')\to E^2E^1(B')\to 0$. 
Hence there are equivalent 
$r$-weights $\omega_B, \omega_{B'}:\Lambda^r\to D$ with kernels $\mbox{Ker}(\omega_B)=B_D\cong B$,  $\mbox{Ker}(\omega_{B'})=B'_D\cong B'$. 
For $\mbox{Ker}(\omega)=B_D$,  the inclusion $B_D\subset \Lambda^r$ induces a 
$\Lambda$-isomorphism $g_B:E^0E^0(B_D)\cong \Lambda^r$. Hence there is a 
$\Lambda$-isomorphism $g_D:E^2E^1(B_D)\to D$ to define an $r$-weight 
$\omega_{B_D}:\Lambda^r\to D$ which is equivalent to $\omega$. 
For equivalent $r$-weights $\omega,  \omega':\Lambda^r\to D$ with 
$\mbox{Ker}(\omega)=B_D$ and $\mbox{Ker}(\omega')=B'_D$,  the five lemma for a 
short exact sequence shows that $B_D$ is $\Lambda$-isomorphic to $B'_D$. 
From finiteness of the $\Lambda$-module $\hom_{\Lambda}(\Lambda^r, D)$,  a desired system of finitely generated torsion-free $\Lambda$-module $B_i\,  (i=1, 2, \dots,  n_r)$ 
of rank $r$ with $E^2E^1(B_i)\cong D$ is obtained. 
Let $\mbox{\boldmath $e$}=\{e_i|\,  i=1, 2, \dots, r\}$ be a standard $\Lambda$-basis of 
$\Lambda^r$. For an $r$-weight $\omega:\Lambda^r\to D$,  
assume that $\omega(e_1)=\omega(e_i)$ for some $i\ne 1$. Then replace the basis element 
$e_i$ with $e_i-e_1$. By continuing this process,  there is a $\Lambda$-isomorphism 
$f_{\Lambda}:\Lambda^r \to \Lambda^r$ such that $\omega'=f_{\Lambda}\omega$ 
is an $r$-weight such that $\omega'$ sends a $\Lambda$-subbase 
$\mbox{\boldmath $e$}'$ of $\mbox{\boldmath $e$}$ injectively and the remaining 
$\Lambda$-subbasis $\mbox{\boldmath $e$}\setminus \mbox{\boldmath $e$}'$ to $0$.
This means that $\mbox{Ker}(\omega)=B_D$ is $\Lambda$-isomorphic to 
$B'\oplus \Lambda^{r-r'}$
for $B'= \Lambda^{r'}\cap \mbox{Ker}(\omega')$ for the $\Lambda$-submodule $\Lambda^{r'}$ given by the $\Lambda$-subbasis $\mbox{\boldmath $e$}'$.
For an $r$-weight $\omega:\Lambda^r\to D$,  
assume that $\omega$ sends $\mbox{\boldmath $e$}$ injectively to $D$. 
Let $n=|D|$. 
Let $\bar \omega:\Lambda^n\to D$ be a $\Lambda$-epimorphism 
extending $\omega$ so that the standard basis 
$\mbox{\boldmath $\bar e$}$ of $\Lambda^n$ bijectively to $D$. 
For every basis element $e_j$ in 
$\mbox{\boldmath $\bar e$}\setminus \mbox{\boldmath $ e$}$,  write
\[\bar x_i(e_j)=\sum_{k=1}^r a_{jk}(t)\xi(e_k)\quad (a_{jk}(t)\in \Lambda).\]
Under the new basis of $\Lambda^n$ obtained by 
replacing every $e_j$ with $e_j-\sum_{k=1}^r a_{jk}(t)e_k$,  the kernel of $\bar \omega:$ is 
$\Lambda$-isomorphic to $B_D\oplus \Lambda^{n-r}$. 
Note that any $n$-weight $\bar \omega:\Lambda^n\to D$ sending the standard basis 
$\mbox{\boldmath $\bar e$}$ of $\Lambda^n$ bijectively to $D$ gives 
the unique class $[\bar \omega]$. Thus,  any $s$-weight $\omega':\Lambda^{s}\to D$ 
sending the standard basis $\mbox{\boldmath $e$}'$ of $\Lambda^n$ injectively to $D$,  
so that there is a $\Lambda$-isomorphism 
$B_D\oplus \Lambda^{n-r}\cong B'_D\oplus \Lambda^{n-r'}$
for $B'_D=\mbox{Ker}(\omega')$.
This completes the proof of Lemma~4.1.

\phantom{x}

\noindent{\bf 4.2: Proof of Theorem~1.1.} 
For the proof of (1),  let $H=R_1(F)$ for a surface-knot $F$ of genus $g$. 
By the second duality of \cite{K86},  
$E^2(R_1(F))=E^2(DA_1(F)/\Theta(F))\cong E^2E^1(BA_2(F)$. 
Hence there is a $t$-anti $\Lambda$-isomorphism 
\[E^2(H)=E^2(R_1(F))\cong E^2E^1(BA_2(F))=E^2E^1(X\oplus Y)\]
by assuming (3). 
Since $E^0E^0(X)$ and $Y$ are free $\Lambda$-modules of rank $g$ and there is a 
$\Lambda$-epimorphism $E^0E^0(X)\to E^2E^1(X)$,  the following inequalities hold. 
\[e(E^2E^1(X\oplus Y))=e(E^2E^1(X))\leq e(E^0E^0(X))=g.\]
Thus,  $e(H)\leq g$,  assuming (3). 
Conversely,  let $H$ be a $(t-1)$-divisible finitely generated $\Lambda$-module with inequality $e(E^2(H))\leq g$. Then $H$ is the first module $A_1(F)$ of a ribbon 
surface-knot $F$ of genus $g$ in $S^4$ with $\Theta(F)=0$ by \cite{K08}. 
Thus, $H=A_1(F)=R_1(F)$,  which shows (1) by assuming (3).
For the proof of (2),  let $[\omega_F]$ be the $g$-class on the finite 
$\Lambda$-module $E^2E^1(BA_2(F))$,  which is $t$-anti $\Lambda$-isomorphic to 
$E^2(R_1(F))$,  so that $[\omega_F]$ is considered as a $g$-class on the finite 
$\Lambda$-module $E^2(R_1(F))$. By Lemmas~4.1,  $BA_2(F)$ is determined by this 
$g$-class on $E^2(R_1(F))$. 
By the first duality of \cite{K86},  the torsion $\Lambda$-module $TA_2(F)=T_DA_2(F)$ 
is $t$-anti $\Lambda$-isomorphic to 
$E^1(T_DH_1(\tilde E, \partial\tilde E;Z))=E^1(T_DA_1(F) )=E^1(T_DR_1(F))=E^1(R_1(F))$ by Lemma~2.1,  
showing (2). 
For the proof of (3),  
note that the zeroth duality of \cite{K86} means that there is a {\it non-singular} 
$\Lambda$-form 
\[S_{\partial}: E^0E^0(BA_2(F))\times E^0E^0BH_2(\tilde E, \partial \tilde E;Z)\to\Lambda\]
extending the non-degenerate $\Lambda$-intersection form 
\[S^B_{\partial}: BA_2(F)\times BH_2(\tilde E, \partial \tilde E;Z)\to\Lambda\]
which also defines a {\it non-degenerate} 
$\Lambda$-Hermitian $\Lambda$-intersection form 
$S^B: BA_2(F)\times BA_2(F)\to\Lambda$. 
Since $DH_2(\tilde M;Z)=0$ and $TH_2(\tilde M;Z)\cong E^1(H_1(\tilde M;Z))$ is 
$(t-1)$-divisible by a method similar to the proof of Lemma~2.1,  
and $H_2(\tilde M, \tilde E;Z)=H_2((V_0, F)\times \mbox{\boldmath $R$};Z)\cong Z^g$,  
there is a natural exact sequence 
\[0\to BA_2(F)\stackrel {i_*}{\to} BH_2(\tilde M;Z)\stackrel {j_*}{\to} 
Z^g\to 0, \]
which induces a natural exact sequence 
\[0\to E^0E^0(BA_2(F))\to E^0E^0(BH_2(\tilde E, \partial \tilde E;Z))\to 
Z^g\to 0.\]
By Lemma~3.1 and Corollary~3.2,  
$BH_2(\tilde M;Z)=X\oplus Y$,  where $Y$ is a free $\Lambda$-module and 
$X$ is characterized by the $\Lambda$-submodule of $BH_2(\tilde M;Z)$ consisting of an element $x$ such that the product 
$u(t)x$ for an element $u(t)\in\Lambda$ with $u(1)=1$ is a linear combination of $x'_i\, (i=1, 2, \dots, g)$. 
If $x\in X$ has $j_*(u(t)x)=0$,  then 
\[j_*(u(t)x)=u(t) j_*(x)=u(1)j_*(x)=j_*(x)=0\]
and $x\in \mbox{Ker}(i_*)=BA_2(F)$. 
Let $X_F=i_*^{-1}(X)$,  which is characterized by the $\Lambda$-submodule of $BA_2(F)$ consisting of an element $x$ such that the product $u(t)x$ for an element $u(t)\in\Lambda$ with $u(1)=1$ is a linear combination of $x'_i=[a(\tilde D'_i)]\, (i=1, 2, \dots, g)$ 
regarded as elements of $BA_2(F)$. This means that $i_*$ defines a $\Lambda$-isomorphism $X_F\cong X$.
Let $Y_F=i_*^{-1}(Y)$ which is a free $\Lambda$-module with basis 
$[a(\tilde D_i)], \,  (i=1, 2, \dots, g)$ 
since $i_*([a(\tilde D_i)])=(t-1)[s(\tilde D_i]$ in $H_2(\tilde M;Z)$. 
This means that $i_*$ defines a natural excat sequence 
\[0\to Y_F\to Y\to Z^g\to 0.\]
Then $BA_2(F)=X_F\oplus Y_F$ 
and the non-degenerate $\Lambda$-Hermitian form 
\[S:E^0E^0(BA_2(F))\times E^0E^0(BA_2(F)) \to \Lambda\]
is given by 
\[S(x_i, x_j)=S(y_i, y_j)=0,  \quad S(x_i, y_j)=(t-1)\delta_{ij} \quad (i, j=1, 2, \dots,  g)\] 
for a $\Lambda$-basis $x_i,  y_i \,  (i =1, 2, \dots,  g)$ of 
$E^0E^0(BA_2(F))=E^0E^0(X_F)\oplus Y_F$ with $x_i \,  (i =1, 2, \dots,  g)$ a $\Lambda$-basis of $E^0E^0(X_F)$ and $y_i\,  (i=1, 2, \dots, g)$ a $\Lambda$-basis of $Y_F$,  
showing (3).
To see (4),  note that $A_3(F)=TA_3(F)$ since $H_3(E;Z)=0$ means that $A_3(F)$ is 
$(t-1)$-divisible. 
By the first duality of \cite{K86},  $T_DA_3(F)$ is $t$-anti $\Lambda$-isomorphic 
to $\hom_{\Lambda}(T_DH_0(\tilde E, \partial \tilde E;Z),  Q(\Lambda)/\Lambda)$ which is $0$.
By the second duality of \cite{K86},  $DA_3(F)$ is $t$-anti $\Lambda$-isomorphic to
$E^1(BH_0(\tilde E, \partial \tilde E;Z))$ which is $0$. Thus,  $A_3(F)=0$,  
showing (4). 
This completes the proof of Theorem~1.1.

\phantom{x}

In the similar way to the proof of (4) in 4.2,  it is shown that 
$H_3(\tilde E, \partial \tilde E;Z)\cong Z$ whose integral generator is  the fundamental class of the infinite cyclic connected covering $\tilde E\to E$  
represented by a leaf of the surface-knot $F$  (see \cite{K98}). 
In fact, by the first duality of \cite{K86},  
$H_3(\tilde E, \partial \tilde E;Z)=T_DH_3(\tilde E, \partial \tilde E;Z)$ which 
is $t$-anti $\Lambda$-isomorphic to $E^1(A_0(F))\cong Z$ . 
The proof of Corollary~1.2 is done as follows. 

\phantom{x}

\noindent{\bf 4.3: Proof of Corollary~1.2.}  Since $\pi$ is  the group of ribbon presentation, 
of deficiency $0$, there is a ribbon torus-knot $T$ in $S^4$ with 
$\pi_1(S^4\setminus T,x_0)=\pi$ and $A_1(T)=D$ (see \cite{K07}).
Since $D$ is a $(t-1)$-divisible finite $\Lambda$-module  with $e(D)=1$, 
the first module $A_1(T_g)$ of $T_g$ in $S^4$ 
is the finite $\Lambda$-module $D^g$,  the direct sum of $g$ copies of $D$,  and 
$E^2(D^g)$ is seen to be $\Lambda$-isomorphic to $D^g$ and $e(D^g)=g$ since $p$ is a prime number. For $p\geq 5$,  the finite $\Lambda$-module $D^g$ 
does not admit any $t$-anti $\Lambda$-automorphism,  so that $\Theta(F)=0$ and 
$A_1(F)=R_1(F)$ for any surface-knot $F$ in $S^4$ with $A_1(F)=D^g$. 
Since $e(R_1(F))=g$, the reduced first module $R_1(F)$ is not $\Lambda$-isomorphic to the reduced first module of any surface-knot of genus $g'<g$ by Theorem~1.1 (1), 
so that $\pi$ is not the fundamental group of any surface-knot of genus $g'<g$. 
 This completes the proof of Corollary~1.2.

\phantom{x}

\noindent{\bf 5. An exact leaf and the torsion-linking of a surface-knot}

Let $V'_F$ be a leaf of a surface-knot $F$ in $S^4$ containing a half surface-basis $D'_i\, (i=1, 2, \dots, g)$ of a surface-basis $D_i, D'_i\, (i=1, 2, \dots, g)$ as proper surfaces. 
Let $W$ be a compact connected oriented 3-manifold with $\partial W=F$,  and 
$V^*=V'_F\cup W$ be the closed oriented 3-manifold obtained from $V_F$ and $W$ by pasting along $F$ with an orientation-reversing diffeomorphism of $F$. 
The following lemma is used for the proof of Theorem~1.4.

\phantom{x}

\noindent{\bf Lemma~5.1.} If $H_1(W;Z)$ is a free abelian group and the loop system $\alpha_i\, (i=1, \dots, g)$ or $\alpha'_i\, (i=1, \dots, g)$ in $F$ represents a basis of the image of the boundary 
homomorphism $\partial_*:H_2(W, F;Z)\to H_1(F;Z)$,  then 
the inclusion $V'_F\to V^*$ induces an isomorphism 
$\mbox{Tor}H_1(V'_F;Z)\to \mbox{Tor}H_1(V^*;Z)$. 

\phantom{x}

\noindent{\bf Proof of Lemmay~5.1.} 
Since the exact leaf $V'_F$ contains the disjoint proper surfaces $C'_i \, (i=1,2,\dots ,g)$, there is a retraction $r_F: V'_F\to \gamma$ for a legged loop system $\gamma$ with the loops 
$\alpha_i\,(i=1,2,\dots,g)$ in $F$ such that the composite $r_F i_F:\gamma\to\gamma$ 
for the inclusion $i_F:\gamma\to V'_F$ is homotopic to the identity. 
Then the homology exact sequence
\[H_2(V'_F,F;Z)\to H_1(\partial F;Z)\to H_1(V'_F;Z)\to H_1(V'_F,F;Z)\to 0\]
induces a split short exact sequence 
\[0\to Z^g\to H_1(V'_F;Z)\to H_1(V'_F,F;Z)\to 0, \]
where $Z^g$ denotes a free abelian group with basis represented by the loops 
$\alpha_i\,(i=1,2,\dots,g)$. Hence 
there are natural isomorphisms 
\[\mbox{Tor}H_1(V'_F;Z)\cong\mbox{Tor}(H_1(V'_F;Z)/\mbox{Im}i_*) \cong 
\cong \mbox{Tor}H_1(V^*, W;Z) \]
for the image $ \mbox{Im}i_*=\mbox{Im}(i_*:H_1(F;Z)\to H_1(V'_F;Z))$. 
Since 
\[H_1(V^*, V'_F;Z)\cong H_1(W, F;Z)\cong H^2(W;Z)\] 
is a free abelian group 
and the image $\mbox{Im}\partial'_*=\mbox{Im}(\partial'_*:H_2(V^*,V'_F;Z) \to 
H_1(V'_F;Z))$ is equal to the image 
$\mbox{Im}i_*\partial_*=\mbox{Im}(i_*\partial_*:H_2(W, F;Z)\to H_1(F;Z)\to H_1(V'_F;Z))$, 
the exact sequence 
\[H_2(V^*, V'_F;Z)\to H_1(V'_F;Z)\to H_1(V^*;Z)\to H_1(V^*, V'_F;Z)\to 0\]
induces a natural isomorphism 
\[\mbox{Tor}(H_1(V'_F;Z)/\mbox{Im}\partial'_*) \cong \mbox{Tor}H_1(V^*;Z).\] 
If the loop system $\alpha_i\, (i=1,2, \dots, g)$ or $\alpha'_i\, (i=1,2, \dots, g)$ represents a basis of 
$\mbox{Im}\partial_*$ in $H_1(F;Z)$,  then there is a natural isomorphism 
\[\mbox{Tor}(H_1(V'_F;Z)/\mbox{Im}\partial_*)\to \mbox{Tor}(H_1(V'_F;Z)/\mbox{Im}i_*).\]
Hence the inclusion $V'_F\to V^*$ induces an isomorphism 
$\mbox{Tor}H_1(V'_F;Z)\to \mbox{Tor}H_1(V^*;Z)$. 
This completes the proof of Lemma~5.1.

\phantom{x}

Theorem~1.4 is shown as follows.

\phantom{x}

\noindent{\bf 5.2: Proof of Theorem~1.4.}
Let $V_F$ be any leaf of $E$ with $\partial V_F=F\times1$ in 
$F\times S^1=\partial E$,  and 
$D_i ,  D'_i\, (i=1, 2, \dots, g)$ any surface-basis of $F$ in $E$ 
with $\partial D_i=\alpha_i\times 1,  \partial D'_i=\alpha'_i\times 1$. 
Let $G_k\, (k=1, 2, \dots, s)$ be closed connected oriented surfaces in $E$ 
representing $\Lambda$-generators of the direct summand $X_F$ of $BA_2(F)$ 
in Theorem~1.1 (2),  which can be disjointedly embedded in $E$ under the covering projection 
$\mbox{proj}: \tilde E\to E$ by \cite[Theorem~4.1]{K00} because $S(X_F, X_F)=0$. 
Since $V_F,  D_i,  G_k$ are all trivially liftable in $\tilde E$,  the leaf $V_F$ is modified so that the interior $\mbox{Int}D_i$ of $D_i$ transversely meets $V_F$ in disjoint simple loops
each of which is null-homologous in $D_i$ and $G_i$ transversely meets $V_F$ in disjoint simple loops each of which is null-homologous in $G_k$. 
Let $D^0_i$ be an innermost piece of the surfaces of $D_i$ divided by the loops 
$V_F\cap D_i$.
Take a normal disk bundle $D^0_i\times D^2$ of $D^0_i$ in $E$ with 
$(\partial D^0_i)\times D^2$ 
a normal disk bundle of the loop $\partial D^0_i$ in $V_F$ and replace 
$(\partial D^0_i)\times D^2$
with $D^0_i\times S^1$ to obtain from $V_F$ to obtain a new leaf of $F$ in $E$.
By continuing this process,  the leaf $V_F$ is modified to have 
$V_F\cap \mbox{Int}D_i=\emptyset\, (i=1, 2, \dots, g)$. 
By continuing the same process after pushing the $\alpha_i\times 1$ into the interior of $V_F$,  
the 3-manifolds $D_i\times S^1\, (i=1, 2, \dots, g)$ are contained in the resulting leaf $V_F$ of $E$. 
By the similar modification,  $V_F$ is modified so that 
a normal circle bundle $G_k\times S^1$ of $G_k$ is made disjoint from $V_F$. 
Replace $V_F$ with a connected sum of $V_F$ and $G_k\times S^1\, (k=1, 2, \dots, s)$ in $E$. 
To show that the resulting leaf $V_F$ is a desired leaf of $F$ in $S^4$,  
let $V=V_F\cup V_0\times 1$ be a closed leaf in the surface-knot manifold 
$M=E\cup V_0\times S^1$,  
where $V_0$ is a handlebody with  a disjoint disk system $d_i\, (i=1, 2, \dots, g)$ 
bounded by the half loop basis $\alpha_i\, (i=1, 2, \dots, g)$.  
Then  the surface $D_i$ extends to a closed surface $s(D_i)=D_i\cup d_i\times 1$ in $V$. 
By 4.2,  $BA_2(F)=X_F\oplus Y_F$,  $BH_2(\tilde M)=X\oplus Y$ and the short 
exact sequence
\[0\to E^0E^0(BA_2(F))\to E^0E^0(BH_2(\tilde M))\to Z^g 
\to 0\]
splits into the isomorphism $X_F\cong X$ and the short exact sequence 
$0\to Y_F\to Y\to Z^{2g} \to 0$. 
Hence 
the natural homomorphism $H_2(\tilde V_F;Z)\to E^0E^0(BH_2(\tilde E))$ with image $X_F$ induces 
the natural homomorphism $H_2(\tilde V;Z)\to E^0E^0(BH_2(\tilde M))$ with image $X$. 
By \cite{K00},  the closed  leaf $V$ of $M$ is a closed exact leaf of $M$,  meaning that the following natural sequence 
\[(*)\qquad 0\to \mbox{Tor}H_2(\tilde M, \tilde V;Z)\to \mbox{Tor}H_1(\tilde V;Z) 
\to \mbox{Tor}H_1(\tilde M;Z)\]
is an exact sequence on integral torsions. 
By Lemma~5.1,  there is a natural isomorphism
$\mbox{Tor}H_1(\tilde V_F;Z) \cong \mbox{Tor}H_1(V;Z)$. 
Since $H_1(\tilde E;Z)$ and $ H_1(\tilde M;Z)$ are $(t-1)$-divisible and 
$(t-1)H_k(\tilde M,\tilde E;Z)=0\, (k=1,2)$,  there is a natural $\Lambda$-isomorphism $H_1(\tilde E;Z) \to H_1(\tilde M;Z)$. 
Further, there is a natural $\Lambda$-isomorphism 
$\mbox{Tor}H_2(\tilde E, \tilde V_F;Z)\to \mbox{Tor}H_2(\tilde M, \tilde V;Z)$.  
In fact, there are $\Lambda$-isomorphisms
$H_2(\tilde E, \tilde V_F;Z)\cong H_2(\tilde E\cup \tilde V, \tilde V;Z)$ and 
$H_k(\tilde M, \tilde E\cup \tilde V;Z)\cong H_k(\tilde V_0\times (I,\partial I);Z)
 \cong H_{k-1}(\tilde V_0;Z)$ by the excision theorem, 
 where $I$ denotes the interval $[0,1]$.  
Since $H_2(\tilde V_0;Z)=0$ and  $H_1(\tilde V_0;Z)\cong\Lambda^g$, 
there is a  natural exact sequence 
\[0\to H_2(\tilde E\cup \tilde V, \tilde V;Z)\to H_2(\tilde M, \tilde V;Z) \to \Lambda^g,\]
which implies a natural $\Lambda$-isomorphism 
 $\mbox{Tor}H_2(\tilde E, \tilde V_F;Z)\cong \mbox{Tor}H_2(\tilde M, \tilde V;Z)$ 
 as desired. Thus, the natural sequence  
\[ 0\to \mbox{Tor}H_2(\tilde E, \tilde V_F;Z)\to \mbox{Tor}H_1(\tilde V_F;Z) 
\to \mbox{Tor}H_1(\tilde E;Z)\]
is equivalent to the exact sequence $(*)$ and $V_F$ is an exact leaf of $E$. 
This completes the proof of Theorem~1.4. 

\phantom{x}

The proof of Corollary~1.5 is given as follows.

\phantom{x}

\noindent{\bf 5.3: Proof of Corollary~1.5.} Since $V$ is a closed exact leaf of $M$,  it is shown in \cite{K00+} that the linking 
$\ell_{V}: \mbox{Tor}H_1(V;Z)\times \mbox{Tor}H_1(V;Z)\to Q/Z$ 
is isomorphic to the orthogonal sum of the torsion linking 
$\ell_M: D^{\theta}H_1(\tilde M;Z) \times D^{\theta}H_1(\tilde M;Z) \to Q/Z$ 
given by the second duality of \cite{K86} and a hyperbolic linking. 
Because $DA_1(F)\cong DH_1(\tilde M;Z)$ and $E^1(BA_2(F))\cong E^1(BH_2(\tilde M;Z))$ 
as $(t-1)$-divisible finite modules, 
the torsion linking $\ell_M$ is $\Lambda$-isomorphic to the torsion linking 
$\ell_F:\Theta(F)\times\Theta(F)\to Q/Z$  by the second duality of \cite{K86}.
By Lemma~5.1, the linking 
$\ell_{V_F}: \mbox{Tor}H_1(V_F;Z)\times \mbox{Tor}H_1(V_F;Z)\to Q/Z$ is non-singular and isomorphic to $\ell_{V}$. Thus,  the linking $\ell_{V_F}$ is an orthogonal sum of $\ell_F$ 
and a hyperbolic linking. 
This completes the proof of Corollary~1.5.

\phantom{x}

\noindent{\bf Acknowledgements}  
This work was partly supported by JSPS KAKENHI Grant Number JP21H00978 and MEXT Promotion of Distinctive Joint Research Center Program JPMXP0723833165.

\end{document}